\def\vers{Sept. 29, 2008, v.5}
\magnification=1200
\hsize=6.5truein
\vsize=8.9truein
\font\bigfont=cmr10 at 14pt
\font\mfont=cmr9
\font\sfont=cmr8
\font\mbfont=cmbx9
\font\sifont=cmti8
\def\scirc{\,\raise.2ex\h{${\scriptstyle\circ}$}\,}
\def\ssb{\raise.2ex\h{${\scriptscriptstyle\bullet}$}}
\def\mopl{\h{$\bigoplus$}}

\def\mcup{\h{$\bigcup$}}
\def\msum{\h{$\sum$}}

\def\b{\beta}
\def\C{{\bf C}}

\def\e{\varepsilon}
\def\g{\gamma}
\def\H{{\bf H}}
\def\HH{\widetilde{H}}
\def\HHH{\widetilde{\bf H}}
\def\cH{{\cal H}}
\def\h{\hbox}
\def\L{{\cal L}}
\def\l{\langle}
\def\N{{\bf N}}

\def\n{\nu}
\def\O{{\cal O}}
\def\P{{\bf P}}
\def\Q{{\bf Q}}
\def\q{\quad}
\def\R{{\bf R}}
\def\r{\rangle}

\def\v{\vee}
\def\X{\widetilde{X}}
\def\XX{{\cal X}}
\def\Y{\widetilde{Y}}
\def\Z{{\bf Z}}
\def\z{\zeta}
\def\zz{\widetilde{\zeta}}
\def\1{{\hskip1pt}}
\def\Coker{{\rm Coker}}
\def\Gr{{\rm Gr}}
\def\Hdg{{\rm Hdg}}
\def\Hom{{\rm Hom}}
\def\HS{{\rm HS}}
\def\MHS{{\rm MHS}}
\def\Ker{{\rm Ker}}
\def\Im{{\rm Im}}
\def\NF{{\rm NF}}
\def\Sing{{\rm Sing}\,}
\def\Spec{{\rm Spec}\,}
\def\ad{{\rm ad}}
\def\can{{\rm can}}

\def\prim{{\rm prim}}
\def\sp{{\rm sp}}

\def\tto{\,\hbox{$\to$}\,}
\def\too{\longrightarrow}

\def\lla{\longleftarrow}

\def\onto{\mathop{\rlap{$\to$}\hskip2pt\hbox{$\to$}}}
\def\simto{\buildrel{\sim}\over\longrightarrow}
\def\({{\rm (}}
\def\){{\rm )}}
\def\1{\hskip1pt}
\h{}
\vskip 1cm
\centerline{\bigfont Generalized Thomas hyperplane sections}

\smallskip
\centerline{\bigfont and relations between vanishing cycles}

\bigskip
\centerline{Morihiko Saito}

\bigskip\medskip
{\narrower\noindent
{\mbfont Abstract.} {\mfont
As is remarked by B.~Totaro, R.~Thomas essentially proved
that the Hodge conjecture is inductively equivalent to the
existence of a hyperplane section, called a generalized
Thomas hyperplane section, such that the restriction to it
of a given primitive Hodge class does not vanish.
We study the relations between the vanishing cycles in the
cohomology of a general fiber, and show that each relation
between the vanishing cycles of type (0,0) with unipotent
monodromy around a singular hyperplane section defines a
primitive Hodge class such that this singular hyperplane
section is a generalized Thomas hyperplane section if and
only if the pairing between a given primitive Hodge class
and some of the constructed primitive Hodge classes does
not vanish. This is a generalization of a construction
by P.~Griffiths.}
\par}

\bigskip\bigskip
\centerline{\bf Introduction}
\footnote{}{{\sifont Date\1}{\sfont:\ \vers}}

\bigskip\noindent
Let $X$ be a smooth complex projective variety of dimension $2n$,
and $\L$ be an ample line bundle on $X$. Let $k$ be a positive
integer such that $\L^k$ is very ample.
Let $S=|\L^k|$, and $\XX$ be the universal family
$\coprod_{s\in S}X_s$ over $S$ with the discriminant $D$.
We assume that the vanishing cycle at a general point of $D$
does not vanish as in [8], XVIII, Cor.~6.4 (replacing $k$ if
necessary).
As is remarked by B.~Totaro (see the last remark in \S 3 of [29]),
R.~Thomas essentially proved that the Hodge conjecture is
inductively equivalent to the existence of a point $0$ of $D$ such
that the restriction $\z|_{X_0}$ of a given primitive Hodge class
$\z$ on $X$ does not vanish (replacing $k$ if necessary).
Here $X_0$ is called a generalized Thomas hyperplane section.
One may assume further that $X_0$ has only ordinary double
points (see loc.~cit.), and $X_0$ is called a Thomas hyperplane
section in this case.
Note that a generalized Thomas hyperplane section is a special
kind of hyperplane section (e.g. it must be reducible if $n=1$).
It has been observed that an explicit construction of a
generalized Thomas hyperplane section for a given primitive
Hodge class is rather difficult.

M.~Green and P.~Griffiths [11] have introduced a notion of
singularities of a normal function.
This is the cohomology class of a normal function.
They showed that non-vanishing of the singularity at $0\in D$ of
the normal function $\n$ associated to $\z$ is equivalent to that
$X_0$ is a Thomas hyperplane section associated to $\z$,
see also [3].
Note that the value $\n_s$ of the normal function at
$s\in S^*:=S\setminus D$ can be viewed as the {\it restriction}
of $\z$ to $X_s$ in the derived category of mixed Hodge
structures using [4].
(This is related to the `restriction' of the Leray spectral
sequence to a fiber in [21], (0.6), see also Remark~(1.2)(i) below.)
Their result shows that the necessary information is not lost
by using this `restriction' even after restricting to a
sufficiently small neighborhood in $S$ of $0\in D$ in the
classical topology.
It implies for example that a Thomas hyperplane section
must have at least two ordinary double points since the
cohomology class of the associated normal function in the
one-variable case is always torsion, see e.g. [22], 2.5.4.
More generally, for a special fiber to be a generalized
Thomas hyperplane section, there must be some relation
between the vanishing cycles in the cohomology of a general
fiber as is shown below.

Let $0\in D$.
To compare the cohomology of $X_0$ with that of $X_s$ for
$s\in S^*:=S\setminus D$, we choose a germ of an irreducible
analytic curve on $S$ whose intersection with $D$ consists of $0$.
Let $C$ be the normalization of the curve.
We assume that $C$ is an open disk.
Let $f:Y\to C$ be the base change of $\XX\to S$ by $C\to S$.
Let $t$ be a local coordinate of $C$ around $0$.
We first assume that $Y_0\,(=X_0)$ has only {\it isolated}
singularities to simplify the exposition.
Then we have the following (see also [26], [27]):

\medskip\noindent
{\bf Proposition~1.} {\it If $\Sing Y_0$ is isolated,
there is an exact sequence of mixed Hodge structures
$$H^{2n-1}(Y_{\infty})\buildrel{\can}\over\too
\mopl_{y\in\Sing Y_0}\,H^{2n-1}(Z_{y,\infty})\too H^{2n}(Y_0)
\buildrel{\sp^{2n}}\over\too H^{2n}(Y_{\infty}),\leqno(0.1)$$
where $H^j(Z_{y,\infty},\Q)$ denotes the vanishing cohomology
at $y\in\Sing Y_0$, and $H^j(Y_{\infty})$ is the cohomology of a
general fiber of $f$ endowed with the limit mixed Hodge structure
at $0\in C$.}

\medskip
Taking the dual of (0.1), we have the dual exact sequence
$$H_{2n-1}(Y_{\infty})\buildrel{\,\,\can^{\v}}\over\lla
\mopl_{y\in\Sing Y_0}\,H_{2n-1}(Z_{y,\infty})\lla
H_{2n}(Y_0)\buildrel{\sp_{2n}}\over\lla H_{2n}(Y_{\infty}).
\leqno(0.2)$$
Set
$$\eqalign{E(Y_0)&=\Ker\big(\sp^{2n}:H^{2n}(Y_0,\Q(n))\to
H^{2n}(Y_{\infty},\Q(n))\big),\cr
R(Y_0)&=\Ker\big(\can^{\v}:\mopl_{y\in\Sing Y_0}\,H_{2n-1}
(Z_{y,\infty},\Q(n))\to H_{2n-1}(Y_{\infty},\Q(n))\big),\cr}$$
where $H_{2n-1}(Y_{\infty},\Q(n))=H^{2n-1}(Y_{\infty},\Q(n))^{\v}$
and similarly for $H_{2n-1}(Z_{y,\infty},\Q(n))$.
Let
$$E^{\v}(Y_0):=E(Y_0)^{\v}=\Coker\big(\sp_{2n}:H_{2n}(Y_{\infty},
\Q(n))\to H_{2n}(Y_0,\Q(n))\big),$$
where $^{\v}$ denotes the dual.
By [3] there is a canonical isomorphism
$$E(Y_0)=\cH^1(j_{!*}\H_{\Q})_0,\leqno(0.3)$$
where $\H$ is a variation of Hodge structure on $S^*:=S\setminus D$
defined by $H^{2n-1}(X_s)(n)$ for $s\in S^*$, and $j_{!*}$ is
the intermediate direct image by the inclusion $j:S^*\to S$,
see [1].
We denote the unipotent monodromy part of $R(Y_0)$ by $R(Y_0)_1$.
For $H=E(Y_0)$, $E^{\v}(Y_0)$, $R(Y_0)_1$, set
$$H^{(0,0)}:=\Hom_{\HS}(\Q,\Gr^W_0H)\,(\supset\Hom_{\MHS}(\Q,H)).$$
This is compatible with the dual.
We say that $E^{\v}(Y_0)^{(0,0)}$ (resp. $E(Y_0)^{(0,0)}$) is the
space of {\it extra Hodge cycles} (resp. {\it cocycles}) on $Y_0$.
An element of $R(Y_0)_1^{(0,0)}$ is called a {\it global relation
between the local vanishing cycles of type $(0,0)$ with unipotent
monodromy} around $Y_0$.
(In the non-islated singularity case, we will omit `local' and
`global'.)

\medskip\noindent
{\bf Theorem~1.}~{\it {\rm (i)}
The restriction of a primitive Hodge class $\z$ to $Y_0$ defines an
extra Hodge cocycle on $Y_0$, i.e. an element of $E(Y_0)^{(0,0)}$.
The latter space is canonically isomorphic to the dual of
$R(Y_0)_1^{(0,0)}$, i.e. there is a canonical isomorphism
$$R(Y_0)_1^{(0,0)}=E^{\v}(Y_0)^{(0,0)}.\leqno(0.4)$$
{\rm (ii)} If $\g_{\b}$ denotes the image of $\b\in R(Y_0)_1^{(0,0)}$
in $H_{2n}(X,\Q(n))^{\prim}$ by the composition of $(0.4)$ with the
canonical morphism
$$E^{\v}(Y_0)\to H_{2n}(X,\Q(n))^{\prim},\leqno(0.5)$$
then $Y_0$ is a generalized Thomas hyperplane section for a
primitive Hodge class $\z$ if and only if
$\langle\z,\g_{\b}\rangle\ne 0$ for some $\b\in R(Y_0)_1^{(0,0)}$.}

\medskip
This is closely related to recent work of M.~Green and
P.~Griffiths~[11].
We are informed that the construction of $\g_{\b}$ was found by
P.~Griffiths ([12], p.~129) in the ordinary double point case,
and the Hodge property of $\g_{\b}$ has been considered by
H.~Clemens (unpublished).
In the general case, using the vanishing cycle functor $\varphi$
in [8], XIII and XIV, we have

\medskip\noindent
{\bf Theorem~2.} {\it Proposition~1 and Theorem~1 hold without
assuming $\Sing Y_0$ is isolated if we replace respectively
$$\mopl_{y\in\Sing Y_0}\,H^{2n-1}(Z_{y,\infty},\Q(n))\q\h{and}\q
\mopl_{y\in\Sing Y_0}\,H_{2n-1}(Z_{y,\infty},\Q(n))$$ by
$$H^{2n-1}(Y_0,\varphi_{f^*t}\Q_Y(n))\q\h{and}\q
H^{2n-1}(Y_0,\varphi_{f^*t}\Q_Y(n))^{\v}.$$}

By (0.3), the dimension $r(Y_0)$ of $R(Y_0)_1^{(0,0)}$ or
$E(Y_0)^{(0,0)}$ is independent of $C$.
So we may assume $C$ smooth for the calculation
of $R(Y_0)_1^{(0,0)}$ and $E(Y_0)^{(0,0)}$, see Remark~(2.8)(i).
As a corollary of Theorem~1, $Y_0$ cannot be a generalized Thomas
hyperplane section if $r(Y_0)=0$.
In the ordinary double point case, the relations are all of
type $(0,0)$ with unipotent monodromy, see Theorem~3 below.
In the isolated singularity case we have a rather explicit
construction of $\g_{\b}$ (which is essentially the same as
Griffiths construction in [12], p.~129), see (2.5) below.
The rank of $\can$ in (0.1) may depend on the
position of the singularities, see Thm.\ (4.5) in [9], p.~208 and
also [10], (3.5).

In the isolated singularity case we have moreover

\medskip\noindent
{\bf Proposition~2.} {\it
If the singularities of $Y_0$ are isolated, then these are
isolated complete intersection singularities,
$\HH^j(Z_{y,\infty})=0$ for $j\ne 2n-1$, and
$\HH^{2n-1}(Z_{y,\infty})$ is independent of $C$ except for the
monodromy.}

\medskip
In the ordinary double point case we show

\medskip\noindent
{\bf Proposition~3.} {\it With the notation of Theorem~1, assume
the singularities of $Y_0$ are ordinary double points.
Then the singularities of the total space $Y$ are of type $A_k$.}

\medskip
Using this, we get the following

\medskip\noindent
{\bf Theorem~3.} {\it With the notation and the assumption of
Proposition~3, the constant sheaf on $Y$ is the intersection
complex up to a shift, i.e.\ $Y$ is a rational homology manifold.
Moreover, the vanishing cohomology at each singular point of $Y_0$
is $\Q(-n)$ as a mixed Hodge structure, and has a unipotent
monodromy.}

\medskip
Combined with [19], Lemma~5.1.4, the first assertion of Theorem~3
implies

\medskip\noindent
{\bf Corollary~1.} {\it With the notation and the assumption of
Proposition~3, let $T$ be the local monodromy around $0$. Then
for $c\in C\setminus\{0\}$ sufficiently near $0$}
$$\Ker\,\can=\Ker(T-id)\q\h{\it on}\,\,\,\, H^{2n-1}(Y_c,\Q).$$

This may be useful in the last section of [3].
Note that Theorem~3 and Corollary~1 do not hold if the
fibers $Y_c$ are even-dimensional with $k$ odd,
see Remark~(2.8)(ii) below.

In Section 1 we review some recent development in the theory of
normal functions, and show certain assertions related to Theorem~1.
In Section 2 we prove the main theorems.

I would like to thank P.~Brosnan, A.~Dimca, P.~Griffiths, J.~Murre
and G.~Pearlstein for useful discussions and valuable comments.
This work is partially supported by Kakenhi 19540023.

\bigskip\bigskip
\centerline{\bf 1. Normal functions}

\bigskip\noindent
{\bf 1.1.~Normal functions associated to primitive Hodge classes.}
With the notation of Introduction let $\H$ be a variation
of Hodge structures of weight $-1$ on $S^*=S\setminus D$ defined
by $H^{2n-1}(X_s,\Z(n))\,(s\in S^*)$.
This gives a family of intermediate Jacobians
$\coprod_{s\in S^*}J^n(X_s)$ containing a constant subfamily
$J^n(X)$.
Take a primitive Hodge class
$$\z\in\Hdg^n(X)^{\prim}\subset H^{2n}(X,\Z(n))^{\prim}.$$
By lifting it to an element of Deligne cohomology and restricting
to $X_s$, it defines an admissible normal function [22]
$$\n\in\NF(S^*,\H)^{\ad}.$$
This is identified with an extension class of $\Z_{S^*}$ by $\H$
as admissible variations of mixed Hodge structures ([15], [28]),
and also with a holomorphic section of $\coprod_{s\in S^*}J^n(X_s)$.
It is well-defined up to a constant section with values in $J^n(X)$.
Let $j:S^*\to S$ denote the inclusion.
The normal function $\n$ has the cohomology class
$$\g(\n)\in H^1(S^*,\H),$$
using the underlying extension class of local systems.
It induces at each $0\in D$
$$\g_0(\n)\in(R^1j_*\H)_0.$$
This is independent of the ambiguity of the normal function.

On the other hand, $\z$ induces by restriction
$$\z|_{X_0}\in H^{2n}(X_0,\Q(n)).$$
Using the functorial morphism $id\to\R j_*j^*$, it induces
further an element of $(R^1j_*\H_{\Q})_0$.
By P.~Brosnan, H.~Fang, Z.~Nie and G.~J.~Pearlstein [3]
(extending the theory of M.~Green and P.~Griffiths [11])
we have the commutativity of the diagram
$$\matrix{\Hdg^n(X)^{\prim}&\too&\NF(S^*,\H)^{\ad}/J^n(X)\cr
\,\,\big\downarrow\raise1.5pt\h{$\scriptstyle\alpha$}&&
{}_{\raise-1.5ex\h{ }}\big\downarrow^{\raise1.5ex\h{ }}\cr
H^{2n}(X_0,\Q(n))&\buildrel{\b}\over\too&(R^1j_*\H_{\Q})_0\cr}
\leqno(1.1.1)$$
and the restriction of $\b$ to the image of $\alpha$ is
injective.

\medskip\noindent
{\bf 1.2.~Remarks.}
(i) The value $\n_s$ of the normal function $\n$ at
$s\in S^*$ may be viewed as the {\it restriction} of a primitive
Hodge class $\z$ to $X_s$ in the derived category of mixed
Hodge structures (using [4]).
The above commutative diagram (1.1.1) asserts that the
restriction of $\z$ to $X_0$ can be calculated by using these
`restrictions' of $\z$ to $X_s$ for $s\in S^*$ sufficiently
near $s$.
This implies that the necessary information is not lost by using
this `restriction' even after restricting to a small neighborhood
of $0$ in the classical topology.
(Note that maximal information will be preserved if we can use
the restriction as algebraic cycles.
This situation is similar to the `restriction' of the Leray
spectral sequence to a fiber in [21], (0.6).)

\medskip
(ii) M.~de Cataldo and L.~Migliorini [6] have proposed a theory
of singularities for primitive Hodge classes using the
decomposition theorem [1] but without normal functions.
For the moment, it is not very clear how to calculate the image of
$\z$ in $(R^1j_*\H_{\Q})_0$ without using the normal functions
as in Remark~(i) above.

\medskip
(iii) A key observation in Thomas argument [29] is that the
algebraic cycle classes coincide with the Hodge classes if and only
if for any Hodge class there is an algebraic cycle class such that
their pairing does not vanish.
For a primitive Hodge class $\z$, the condition that
$\z|_{X_0}\ne 0$ for some $0\in D$ implies the existence of an
algebraic cycle such that their pairing does not vanish as in
Remark~(iv) below.
Note, however, that this condition does not immediately imply the
algebraicity of $\z$ (unless it is satisfied for any $\z$)
since this is insufficient to show the coincidence of the
algebraic and Hodge classes.

\medskip
(iv) As is remarked by B.~Totaro (see the last remark in \S 3 of
[29]), Thomas argument is extended to the case of arbitrary
singularities by using the injectivity of
$$\Gr^W_{2n}H^{2n}(X_0,\Q)\to H^{2n}(\widetilde{X_0},\Q),$$
where $\widetilde{X_0}\to X_0$ is a desingularization.
(This injectivity follows from the construction of mixed Hodge
structure using a simplicial resolution [7]).
If $\z|_{X_0}\ne 0$ for a primitive Hodge class $\z$,
then we can apply the Hodge conjecture for $\widetilde{X_0}$ as an
inductive hypothesis to construct an algebraic cycle on $X_0$
whose pairing with $\z$ does not vanish,
using the above injectivity (together with the strict
compatibility of the weight filtration $W$).
This point is the only difference between the general case and
the ordinary double point case in Thomas argument~[29], and
the hypothesis on ordinary double points is not used
in the other places (as far as the proof of the Hodge conjecture
is concerned).

\medskip\noindent
{\bf 1.3.~Cohomology classes of normal functions.}
Let $S$ be a complex manifold, and $S^*$ be an open subset such
that $D:=S\setminus S^*$ is a divisor.
Let $\H$ be a polarizable variation of Hodge structure of weight
$-1$ on $S^*$. Let
$$\n\in\NF(S^*,\H_{\Q})_S^{\ad}:=
\NF(S^*,\H)_S^{\ad}\otimes_{\Z}\Q.$$
It is an extension class of $\Q$ by $\H_{\Q}$ as admissible
variations of mixed $\Q$-Hodge structures ([15], [28]),
and is identified with an extension class as shifted mixed Hodge
modules on $S$
$$\Q_S\to\R j_*\H_{\Q}[1].\leqno(1.3.1)$$
Here $\Q_S$ and $\R j_*\H_{\Q}$ are mixed Hodge modules up to a shift
of complex by $n$ since $D$ is a divisor.
Let $j_{!*}\H_{\Q}$ be the intermediate direct image, i.e.\
the intersection complex up to a shift of complex by $n$, see [1].
Then (1.3.1) factors through $(j_{!*}\H_{\Q})[1]$ by the
semisimplicity of the graded pieces of mixed Hodge modules
since the weight of $\H$ is $-1$, see [3], [18].

Let $i_0:\{0\}\to S$ denote the inclusion.
Then (1.3.1) induces a morphism of mixed Hodge structures
$$\Q\to H^1i_0^*\R j_*\H_{\Q},$$
factorizing through $H^1i_0^*j_{!*}\H_{\Q}$.
The image of $1\in\Q$ by this morphism is called the
cohomology class of $\n$ at $0$.
We get thus the morphisms
$$\NF(S^*,\H_{\Q})_S^{\ad}\to\Hom_{\MHS}(\Q,H^1i_0^*j_{!*}\H_{\Q})
\hookrightarrow\Hom_{\MHS}(\Q,H^1i_0^*\R j_*\H_{\Q}).
\leqno(1.3.2)$$
Here the injectivity of the last morphism easily follows from the
support condition on the intersection complexes, see [3] (and also
(1.4) below for the normal crossing case).

\medskip\noindent
{\bf 1.4.~Intersection complexes in the normal crossing case.}
With the above notation, assume that $S$ is a polydisk $\Delta^n$
with coordinates $t_1,\dots,t_n$, $S^*=(\Delta^*)^n$,
and the local monodromies $T_i$ around $t_i=0$ are unipotent.
Let $H$ be the limit mixed Hodge structure of $\H$, see [24].
Set $N_i=\log T_i$.
The functor $i_0^*$ between the derived category of mixed Hodge
modules [20] is defined in this case by iterating the mapping
cones of
$$\can:\psi_{t_i}\to\varphi_{t_i}.$$
So $H^1i_0^*\R j_*\H_{\Q}$ is calculated by the cohomology at
degree 1 of the Koszul complex
$$K^{\ssb}(H;N_1,\dots,N_n):=\big[0\to H\buildrel{\oplus_iN_i}\over
{\too}\mopl_i\,H(-1)\to\mopl_{i\ne j}\,H(-2)\to\cdots\big],$$
where $H$ is put at the degree 0.
Moreover, it is known (see e.g.\ [5]) that $H^1i_0^*j_{!*}\H_{\Q}$
is calculated by the cohomology at degree 1 of the subcomplex
$$I^{\ssb}(H;N_1,\dots,N_n):=\big[0\to H\buildrel{\oplus_iN_i}\over
{\too}\mopl_i\,\Im\,N_i\to\mopl_{i\ne j}\,\Im\,N_iN_j\to\cdots\big].$$
Define $(\mopl_i\,\Im\,N_i)^0=\Ker(\mopl_i\,\Im\,N_i\to
\mopl_{i\ne j}\,\Im\,N_iN_j)$ so that
$$(\mopl_i\,\Im\,N_i)^0/\Im(\mopl_i\,N_i)=
H^1I^{\ssb}(H;N_1,\dots,N_n)=H^1i_0^*j_{!*}\H_{\Q}.
\leqno(1.4.1)$$

\medskip\noindent
{\bf 1.5.~Remark.} With the notation and the assumption
of $(1.4)$, assume $\H$ is a nilpotent orbit.
Then it is easy to show (see e.g.\ [23]) that $(1.3.2)$ induces a
surjective morphism
$$\NF(S^*,\H_{\Q})_S^{\ad}\onto\Hom_{\MHS}(\Q,
H^1i_0^*j_{!*}\H_{\Q}).\leqno(1.5.1)$$

\medskip\noindent
{\bf 1.6.~Remark.} If $\H$ is not a nilpotent orbit, let
$\HHH$ denote the associated nilpotent orbit.
The target of (1.5.1) does not change by replacing $\H$ with
$\HHH$.
However, the image of (1.5.1) can change in general
(see [23]).

\medskip
The following is closely related to Theorem~1 in the case where
$X_0$ has only ordinary double points and $D$ is a divisor with
normal crossings around $0\in S$ (since $\Im\,N_i$ is generated
by a vanishing cycle via the Picard-Lefschetz formula,
see [8], XV, Th.\ 3.4).
In the geometric case, this is due to [11].

\medskip\noindent
{\bf 1.7.~Proposition.} {\it With the notation and the assumption
of $(1.4)$, assume $N_iN_j=0$ for any $i,j$, and
$\Im\,N_i\subset H(-1)$ is a direct sum of $1$-dimensional
mixed Hodge structures for any $i$.
Let $r$ be the dimension of the relations between the $\Im\,N_i$,
i.e.
$$r=\dim\Ker\big(\mopl_{i=1}^n\,\Im\,N_i\to H(-1)\big).
\leqno(1.7.1)$$
Then}
$$\dim\Hom_{\MHS}(\Q,(R^1j_*\H_{\Q})_0)=
\dim\cH^1(j_{!*}\H_{\Q})_0=r.\leqno(1.7.2)$$

\medskip\noindent
{\it Proof.}
Since $N_i^2=0$, the weight filtration of $\Im\,N_i\subset H(-1)$
is given by the monodromy filtration for $\sum_{j\ne i}N_j$,
see [5] and the references there.
Then $\Im\,N_i$ has type $(0,0)$ (using the hypotheses on $N_i$),
and the first isomorphism of (1.7.2) follows.

There is a nondegenerate pairing $\l *,*\r$ of $H$ giving a
polarization of mixed Hodge structure; in particular
$\l N_iu,v\r=-\l u,N_iv\r$.
It is well-known (see e.g.\ loc.~cit.) that there is a
nondegenerate pairing $\l *,*\r_i$ of $\Im\,N_i$ defined by
$$\l N_iu,N_iv\r_i:=\l N_iu,v\r=-\l u,N_iv\r.$$
Then the morphisms
$$\mopl_i\,N_i:H\to\mopl_i\,\Im\,N_i,\q
\mopl_i\,\Im\,N_i\to H(-1)$$
are identified with the dual of each other, and
$\mopl_i\,\Im\,N_i=(\mopl_i\,\Im\,N_i)^0$ in the notation of
(1.4.1).
So the assertion follows.

\vfill\eject
\centerline{\bf 2. Vanishing cycles}

\bigskip\noindent
{\bf 2.1.~Proof of Proposition~1 in the general case.}
We show Proposition~1 without assuming $\Sing Y_0$ is isolated
as in Theorem~2.
Forgetting the mixed Hodge structure, this is more or less
well-known, see [8], XIII and XIV.
For the compatibility with the mixed Hodge structure,
we can argue as follows.
(If $\Sing Y_0$ is isolated, we can use [26], [27].)

Since $f$ is projective and $C$ can be replaced by a sufficiently
small open disk, we may assume that $Y$ is an intersection of
divisors on $\P^m\times C$.
Then $\Q_Y$ is defined in the derived category of mixed Hodge
modules, see e.g. the proof of Cor.~2.20 in [20].
(In this case, $Y$ is a complete intersection and $\Q_Y[2n]$ is
a perverse sheaf so that it underlies a mixed Hodge module.)
Let $t$ be a local coordinate around $0\in C$, and $i:Y_0\to Y$
be the inclusion.
Then there is a distinguished triangle in the derived categories
of mixed Hodge modules on $Y_0$
$$i^*\Q_Y\too\psi_{f^*t}\Q_Y\too\varphi_{f^*t}\Q_Y
\buildrel{+1}\over\too.$$
Taking the direct image of this triangle by the morphism
$Y_0\to pt$, the assertion follows.

\medskip\noindent
{\bf 2.2.~Proof of Proposition~2.}
This follows from the theory of versal flat deformations
of complete intersections with isolated singularities in the
category of analytic spaces (see [14], [30]) using the base
change of Milnor fibrations.
(The vanishing for $j\ne 2n-1$ follows also from the fact that
$\Q_Y[2n]$ and $\varphi_{f^*t}\Q_Y[2n-1]$ are perverse sheaves
since $Y$ is a complete intersection.)

For each singular point $y_i$, we see that $(Y_0,y_i)$ is
a complete intersection since $\XX$ is smooth, and hence there is
a versal flat deformation of $(Y_0,y_i)$
$$h_i:(\C^{n_i},0)\to(\C^{m_i},0),\leqno(2.2.1)$$
such that $(Y,y_i)\to(C,0)$ is isomorphic to the base change of
$h_i$ by a morphism
$$\rho_i:(C,0)\to(\C^{m_i},0).$$

Let $B_i,B'_i$ be open balls in $\C^{n_i},\C^{m_i}$ with radius
$\e_i$ and $\e'_i$ respectively.
Let $D'_i\subset B'_i$ be the discriminant of $h_i$.
For $1\gg\e_i\gg\e'_i>0$, consider the restriction of $h_i$
$$B_i\cap h_i^{-1}(B'_i\setminus D'_i)\to B'_i\setminus D'_i.$$
This is a $C^{\infty}$ fibration, and the fiber
$B_i\cap h_i^{-1}(s)$ for $s\in B'_i\setminus D'_i$ is
topologically independent of $1\gg\e_i\gg\e'_i>0$.
We have moreover for $s\in B'_i\setminus D'_i$ (see [13], [16])
$$\HH^j(B_i\cap h_i^{-1}(s),\Q)=0\q\h{for}\,\,\,j\ne 2n-1.$$
Using the base change of this fibration by $\rho_i$, the assertion
follows.

\medskip\noindent
{\bf 2.3.~Proof of Proposition~3.}
This follows from the theory of versal flat deformations
explained in (2.2).
Indeed, by the assumption that the singularities of $Y_0$ are
ordinary double points, we have $m_i=1$ and $h_i$ in (2.2.1) is
given by
$$h:(\C^{2n},0)\ni(x_1,\dots,x_{2n})\mapsto\msum_{i=1}^{2n}x_i^2
\in(\C,0).\leqno(2.3.1)$$
If the degree of $\rho_i:(C,0)\to(\C,0)$ is $k_i+1$ with $k_i\in\N$,
then $(Y,y_i)$ is locally isomorphic to a hypersurface defined by
$$\msum_{i=1}^{2n}\,x_i^2=t^{k_i+1},$$
where $t$ is a local coordinate of $C$.
So it has a singularity of type $A_{k_i}$ if it is singular.

\medskip\noindent
{\bf 2.4.~Proof of Theorem~1 in the general case.}
We show Theorem~1 in the general case as in Theorem~2.
The first assertion follows from the hypothesis that $\z$ is
Hodge and primitive.
By (0.2) modified as in Theorem~2,
$R(Y_0)_1^{(0,0)}$ is canonically isomorphic to
$E^{\v}(Y_0)^{(0,0)}$, and this is the dual of $E(Y_0)^{(0,0)}$.
Thus Theorem~1\1(i) is proved in the general case.

For $\b\in R(Y_0)_1^{(0,0)}$, let $\g'_{\b}$ be the corresponding
element in $E^{\v}(Y_0)^{(0,0)}$.
We have the canonical morphism (0.5) using the Lefschetz
decomposition for $X$ since the image of $H_{2n}(Y_{\infty},\Q(n))$
is contained in the non-primitive part.
We define $\g_{\b}$ to be the image of $\g'_{\b}$ by (0.5).
Here the pairing with $\z$ does not change by taking only
the primitive part, since the pairing between the primitive
part and the non-primitive part vanishes.
So we get Theorem~1\1(ii) in the general case using Theorem~1\1(i).

\medskip\noindent
{\bf 2.5.~Construction of $\g_{\b}$ in the isolated singularity
case.}
We can construct $\g_{\b}$ in Theorem~1 rather explicitly in this
case as follows (forgetting the mixed Hodge structure).
For $c\in C$ sufficiently near $0\in C$, let $\rho:Y_c\to Y_0$ be
a good retraction inducing an isomorphism over
$Y_0\setminus\Sing Y_0$.
(This can be constructed by taking an embedded resolution and
composing it with a good retraction for the resolution,
see also [8], XIV.) Set
$$Z_c=\mcup_{y\in\Sing Y}\,\rho^{-1}(y)\cap Y_c.$$
Since $H^j(Y_c,Z_c)=H_c^j(Y_c\setminus Z_c)$, there are
isomorphisms
$$\rho^*:H^j(Y_0,Z_0)\simto H^j(Y_c,Z_c)\q\h{for any}
\,\,j,$$
and $H^j(Y_0,Z_0)=H^j(Y_0)$ for $j\ge 2$.
So the exact sequence (0.1) is identified with
$$H^{2n-1}(Y_c)\to H^{2n-1}(Z_c)\to H^{2n}(Y_c,Z_c)\to
H^{2n}(Y_c),$$
and similarly for the dual.
Take a topological relative cycle $\g'\in H_{2n}(Y_c,Z_c)$
whose image in $H_{2n-1}(Z_c)$ is $\b$. Then
$$\rho_*\g'\in H_{2n}(Y_0,Z_0)=H_{2n}(Y_0),$$
and $\g_{\b}$ in Theorem~1 is the primitive part
of its image in $H_{2n}(X)$.
This construction is essentially the same as the one found by
P.~Griffiths ([12], p.~129) in the ordinary double point case.

\medskip\noindent
{\bf 2.6.~Remark.}
In case $n=1$, the above construction is quite intuitive since we
get a topological 2-chain bounded by vanishing cycles on a nearby
fiber $Y_c$, which gives an algebraic cycle supported on the
singular fiber $Y_0$ by taking the direct image by $\rho$.
However, this does not immediately imply the Hodge conjecture
for this case since the problem seems to be converted to the
one studied in [29] using the pairing between Hodge classes
and algebraic cycles.
The situation may be similar for $n\ge 2$ if one assumes the Hodge
conjecture for a desingularization of $Y_0$.

\medskip\noindent
{\bf 2.7.~Proof of Theorem~3.}
A hypersurface singularity is a rational homology manifold
if and only if $1$ is not an eigenvalue of the Milnor monodromy.
In the isolated singularity case this follows from the Wang
sequence, see e.g. [17].
It is also well-known (see loc.~cit.) that the eigenvalues of the
Milnor monodromy of an even-dimensional $A_k$-singularity are
$$\exp(2\pi ip/(k+1))\q\h{with}\q p=1,\dots,k.$$
(This is a simple case of the Thom-Sebastiani formula [25].)
So the first assertion follows.

For the last assertion, recall that the weight filtration
on the unipotent monodromy part of $\varphi_{f^*t}\Q_Y[2n-1]$ is
the monodromy filtration shifted by $2n$ so that the middle
graded piece has weight $2n$, see [20].
Using the base change of the Milnor fibration by $\rho_i$, we see
that the vanishing cohomology is 1-dimensional and has a
unipotent monodromy in this case.
So the vanishing cohomology is pure of weight $2n$,
and the assertion follows.

\medskip\noindent
{\bf 2.8.~Remarks.}
(i) In the isolated singularity case, we can choose a curve
$C\subset S$ passing through $0$ and such that the base change $Y$
of $\XX$ is smooth by using a linear system spanned by $X_0$ and
$X_s$ such that $X_s$ does not meet any singular points
of $X_0$ (as is well-known).
In this case Proposition~1 follows from the theory of
Steenbrink [26].
However, it is sometimes desirable to show Proposition~1 for
$C\subset S$ such that $Y$ is not smooth, e.g.\ when
$C$ is the image of a curve on a resolution of singularities
of $(S,D)$, see the last section of [3].

\medskip
(ii) Theorem~3 and Corollary~1 do not hold if the fibers are
$2n$-dimensional and if the singularities are of type $A_k$
with $k$ odd.
In this case $\Q_Y[\dim Y]$ is not an intersection complex,
and the monodromy $T$ on $H^{2n}(Y_c,\Q)$ is the identity
since the $k$ are odd.
However, we have non-vanishing of the canonical morphism
$$\can:H^{2n}(Y_c,\Q)\to\mopl_{y\in\Sing Y_0}\Q(-n),$$
for example, if it is obtained by the base change under a double
covering $C\to C'$ of a morphism $Y'\to C'$ with $Y'$ smooth.

\medskip
(iii) It is known that the rank of the morphism $\can$ in
Proposition~1 may depend on the position of the singularities,
see e.g.\ Thm. (4.5) in [9], p.~208 and also [10], (3.5).
Here the examples are hypersurfaces in $\P^{2n}$.
One can construct a hypersurface $X$ in $\P^{2n+1}$ whose hyperplane
section is a given hypersurface $Y$ as follows.

Let $f$ be an equation of $Y$, which is a homogeneous polynomial
of degree $d$.
Let $g=\msum_{i=0}^d\,g_i$, where $g_i$ is a homogeneous polynomial
of degree $i$, and $g_d=f$.
Let $X$ be the closure of $\{g=0\}\subset\C^{2n+1}$ in $\P^{2n+1}$.
Then $X$ is smooth along its intersection with the divisor at
infinity $\P^{2n}$ if $\{g_{d-1}=0\}$ does not meet the
singularities of $Y=\{g_d=0\}$.
As for the intersection of $X$ with the affine space $\C^{2n+1}$,
it is defined by $g$, and is smooth if $g_0$ is sufficiently
general since the critical values of $g$ are finite.
(It does not seem easy to construct $X$ having two given hyperplane
sections.
If we consider a pencil defined by a linear system spanned by two
hypersurfaces we get a pencil of the projective space embedded by
$\O(d)$ in a projective space.)

\medskip
Using Remark~(2.8)(i) above we can show the following
(which would be known to specialists).

\medskip\noindent
{\bf 2.9.~Proposition.} {\it For an ordinary double point $x$
of $X_0$, let $(\Sigma,x)$ be the critical locus near $x$, and
$(D_x,0)$ be its image in $S$.
Then $(\Sigma,x)$ is isomorphic to $(D_x,0)$ and they are smooth.}

\medskip\noindent
{\it Proof.}
By [14], [30], there is a morphism
$$g_x:(S,s)\to(\C,0),$$
such that $(\XX,x)\to(S,s)$ is isomorphic to the base change of
$h_i$ in (2.2.1) by $g_x$.
Then we have $D_x=g_x^{-1}(0)$.
Let $i:(C,0)\to(S,s)$ be a curve in Remark~(2.8)(i).
The composition $g_x\scirc i$ has degree 1 since the base change
of $h$ by it has otherwise a singularity.
So $g_x$ has a section and hence $g_x$ and $D_x$ are smooth.
Then $(\Sigma,x)$ is also smooth since $(\Sigma,x)\to(D_x,s)$ is
bijective.
Thus the assertion is proved.

\medskip\noindent
{\bf 2.10.~Remarks.}
(i) Let $x_1,\dots,x_m$ be ordinary double points on $X_0$.
Then we have a morphism
$$G:(S,0)\to(\C^m,0),$$
whose composition with the $i$-th projection $pr_i:\C^m\to\C$
coincides with $g_{x_i}$ in the proof of Proposition~(2.9).
It is not easy to calculate $G$ although $g_{x_i}=pr_i\scirc G$
is smooth by Proposition~(2.9).

\medskip
(ii) Assume $X_0=Y_0$ has only ordinary double points as
singularities, and let $\Y_0\to Y_0$ be the resolution of
singularities obtained by blowing up along all the singular points
$x_i\,\,(i=1,\dots,m)$.
In this case we have by [11]
$$\eqalign{r(Y_0)&=h^{n,n-1}(\Y_0)-h^{n,n-1}(Y_{\infty})+m,\cr
&=\dim\Hdg^{n-1}(\Y_0)-\dim\Hdg^{n-1}(X)+(1-\delta_{n,1})m.}$$
Since $H^{2n-2}(Y_0)=H^{2n-2}(Y_{\infty})=H^{2n-2}(X)$ and
$H^{2n-2}(E_i)(n-1)\cong\Q^2$, they are closely related to the
exact sequences
$$\eqalign{&0\tto H^{2n-2}(Y_0)\tto H^{2n-2}(\Y_0)\tto\mopl_{i=1}^m
\,\HH^{2n-2}(E_i)\to H^{2n-1}(Y_0)\tto H^{2n-1}(\Y_0)\tto 0,\cr
&0\tto H^{2n-1}(Y_0)\tto H^{2n-1}(Y_{\infty})\tto\mopl_{i=1}^m\,
\Q(-n)\tto E(Y_0)(-n)\tto 0,}$$
since these two imply also $h^{n,n-1}(Y_0)=h^{n,n-1}(\Y_0)$ and
$$\dim\Gr^W_{2n-2}H^{2n-1}(Y_0)=\dim\Gr^W_{2n}H^{2n-1}(Y_{\infty})
=m-r(Y_0).$$

\medskip\noindent
{\bf 2.11.~Remarks.}
(i) Let $\X\to\P^1$ be a Lefschetz pencil where
$\pi:\X\to X$ is the blow-up along the intersection of two general
hyperplane sections.
Let $X_t$ be a general fiber with the inclusion $i_t:X_t\to\X$.
If $2p<\dim X$, then the Leray spectral sequence for the Lefschetz
pencil induces an exact sequence
$$0\to H^{2p-2}(X_t,\Q)(-1)\buildrel{(i_t)_*}\over\too
H^{2p}(\X,\Q)\buildrel{i_t^*}\over\too H^{2p}(X_t,\Q).
\leqno(2.11.1)$$
This can be used to solve a minor problem in an argument in [29].
Indeed, by a Hilbert scheme argument (using the countability
of the irreducible components of the Hilbert scheme),
one can construct an algebraic cycle class $\xi$ with rational
coefficients on $\X$ whose restriction to $X_t$ coincides with
the restriction to $X_t$ of a given primitive Hodge class $\z$
on $X$ where $t\in\P^1$ is quite general.
However, it is not very clear whether $\xi=\pi^*\z$ in loc.~cit.
This problem can be solved by considering the difference
$\pi^*\z-\xi$ since it is a Hodge class and belongs to the image
of $(i_t)_*$ by (2.11.1) so that the inductive hypothesis on the
Hodge conjecture applies.
(This argument seems to be simpler than the one given by
M.~de Cataldo and L.~Migliorini [6].)

\medskip
(ii) The Hilbert scheme argument in [29] can be replaced by
`spreading out' of cycles (a technique initiated probably by
S.~Bloch [2], see also [31]).
Indeed, let $k$ be an algebraically closed subfield of $\C$
which has finite transcendence degree and over which
the Lefschetz pencil $\X_k\to\P^1_k$ is defined.
Let $U$ be a dense open subvariety of $\P_k^1$ over which the
fibers are smooth.
Let $t$ be a $k$-generic point of $\P^1_{\C}$.
Using the inductive hypothesis, the restriction of a Hodge cycle
$\z$ to $X_t$ is represented by an algebraic cycle with rational
coefficients $\xi_t$.
This $\xi_t$ is defined over a subfield $K$ of $\C$ which contains
$k(t)$ and is finitely generated over $k$.
Let $R$ be a finitely generated $k$-subalgebra of $K$ whose
quotient field is $K$ and such that $\xi_t$ is defined over $R$.
Let $\X_{k,V}$ denote the base change of $\X_k\to\P^1_k$ by
$V:=\Spec R\to\P_k^1$,
where we may assume that $V\to\P_k^1$ factors through $U$.
Then $\xi_t$ is defined on $\X_{k,V}$, and its cycle class is
defined as a global section of the local system on $V_{\C}$,
and coincides with the pull-back of the global section
$\zz$ on $U_{\C}\subset\P^1_{\C}$ which is defined by
the restrictions $\z|_{X_{t'}}$ for $t'\in U_{\C}$.
(Indeed, $V$ and $V_{\C}$ are irreducible, and the two global
sections on $V_{\C}$ coincide at the point of $V_{\C}$ determined
by the inclusion $R\to\C$).
Taking a curve $C$ on $V$ which is dominant over $\P^1_k$,
and using the direct image by the base change of $C\to U$
(and dividing it by the degree of $C\to U$), we get a cycle
on $\X_{k,U}\subset\X_k$ whose cycle class coincides with
$\zz$ as global sections on $U_{\C}$, where we may
assume that $C$ is finite over $U$ replacing $U$ and $C$ if
necessary.
Then we can extend it to a cycle on $\X_k$ by taking the closure.

\medskip
(iii) The above argument is essentially explained in
Remarks~(1.3)(ii) and (1.10)(ii) of [21], where it is noted
that if HC$(X,p)$ denotes the Hodge conjecture for cycles of
codimension $p$ on a smooth projective variety $X$,
then HC$(X,p)$ for $p>\dim X/2$ is reduced to HC$(Y,p-1)$ for
a smooth hyperplane section $Y$ (using the Gysin morphism together
with the weak Lefschetz theorem), and for $p<\dim X/2$, it is reduced
to HC$(Y,p)$ and HC$(Y,p-1)$ for a quite general hyperplane section
$Y$ (using a Lefschetz pencil $\X\to\P^1$ and spreading out as above).
Moreover, the problem in Remark~(2.11)(i) above is also mentioned at
the end of Remark~(1.3)(ii) in loc.~cit.\
(i.e.\ HC$(Y,p-1)$ is necessary in the second case).

\bigskip\bigskip
\centerline{{\bf References}}

\medskip
{\mfont
\item{[1]}
A.~Beilinson, J.~Bernstein and P.~Deligne, Faisceaux pervers,
Ast\'erisque, vol. 100, Soc. Math. France, Paris, 1982.

\item{[2]}
S.~Bloch, Lectures on algebraic cycles, Duke University
Mathematical series 4, Durham, 1980.

\item{[3]}
P.~Brosnan, H.~Fang, Z.~Nie and G.~J.~Pearlstein,
Singularities of admissible normal functions
(arXiv:0711.0964).

\item{[4]}
J.~Carlson, Extensions of mixed Hodge structures, in Journ\'ees
de G\'eom\'etrie Alg\'ebrique d'Angers 1979, Sijthoff-Noordhoff
Alphen a/d Rijn, 1980, pp. 107--128.

\item{[5]}
E.~Cattani, A.~Kaplan, W.~Schmid, $L^2$ and intersection
cohomologies for a polarizable variation of Hodge structure,
Inv.\ Math.\ 87 (1987), 212--252.

\item{[6]}
M.~de Cataldo and L.~Migliorini, A remark on singularities of
primitive cohomology classes (arXiv:0711.1307).

\item{[7]}
P.~Deligne, Th\'eorie de Hodge, III, Publ.\ Math.\ IHES 44
(1974), 5--77.

\item{[8]}
P.~Deligne and N.~Katz, Groupes de Monodromie en G\'eom\'etrie
Alg\'ebrique, S\'eminaire de G\'eom\'etrie Alg\'ebrique du
Bois-Marie 1967--1969 (SGA~7~II), Lect.\ Notes in Math.,
vol.~340, Springer, Berlin, 1973.

\item{[9]}
A.~Dimca, Singularities and Topology of Hypersurfaces, Springer,
Berlin 1992.

\item{[10]}
A.~Dimca, M.~Saito and L.~Wotzlaw, A generalization of Griffiths
theorem on rational integrals, II (preprint, arXiv:math/0702105),
to appear in Michigan Math.\ J.

\item{[11]}
M.~Green and P.~Griffiths, Algebraic cycles and singularities
of normal functions, in Algebraic Cycles and Motives, Vol.~1
(J.~Nagel and C.~Peters, eds) London Math.\ Soc.\ Lect.\ Note
Series 343, Cambridge University Press 2007, 206--263.

\item{[12]}
P.~Griffiths, A theorem concerning the differential equations
satisfied by normal functions associated to algebraic cycles,
Amer.\ J.\ Math.\ 101 (1979), 94--131.

\item{[13]}
H.~Hamm, Lokale topologische Eigenschaften komplexer R\"aume,
Math.\ Ann.\ 191 (1971), 235--252.

\item{[14]}
A.~Kas and M.~Schlessinger, On the versal deformation of a complex
space with an isolated singularity, Math.\ Ann.\ 196 (1972), 23--29.

\item{[15]}
M.~Kashiwara, A study of variation of mixed Hodge structure,
Publ.\ RIMS,\ Kyoto Univ. 22 (1986), 991--1024.

\item{[16]}
E.~Looijenga, Isolated singular points on complete intersections,
London Math.\ Soc.\ Lect.\ Note Ser.\ 77,
Cambridge University Press, 1984.

\item{[17]}
J.~Milnor, Singular points of complex hypersurfaces, Annals of
Mathematics Studies No. 61, Princeton University Press 1968.

\item{[18]}
C.~Peters and M.~Saito, Lowest weights in cohomology of variations
of Hodge structure (preprint).

\item{[19]}
M.~Saito, Modules de Hodge polarisables, Publ. RIMS, Kyoto Univ.\
24 (1999), 849--995.

\item{[20]}
M.~Saito, Mixed Hodge modules, Publ. RIMS, Kyoto Univ.\ 26
(1990), 221--333.

\item{[21]}
M.~Saito, Some remarks on the Hodge type conjecture, in Motives
(Seattle, WA, 1991), Proc.\ Sympos.\ Pure Math., 55, Part 1,
Amer.\ Math.\ Soc., Providence, RI, 1994, pp.~85--100.

\item{[22]}
M.~Saito, Admissible normal functions, J. Alg.\ Geom.\ 5 (1996),
235--276.

\item{[23]}
M.~Saito, Variations of Hodge structures having no cohomologically
nontrivial admissible normal functions, preprint.

\item{[24]}
W.~Schmid, Variation of Hodge structure: The singularities of the
period mapping, Inv. Math. 22 (1973), 211--319.

\item{[25]}
M.~Sebastiani and R.~Thom, Un r\'esultat sur la monodromie,
Inv.\ Math.\ 13 (1971), 90--96.

\item{[26]}
J.H.M.~Steenbrink, Mixed Hodge structure on the vanishing
cohomology, in Real and Complex Singularities (Proc. Nordic
Summer School, Oslo, 1976) Alphen a/d Rijn: Sijthoff \& Noordhoff
1977, pp. 525--563.

\item{[27]}
J.H.M.~Steenbrink, Mixed Hodge structures associated with isolated
singularities, in Proc.\ Symp.\ Pure Math.\ 40, Part 2,
Amer.\ Math.\ Soc., Providence, RI, 1983, pp.~513--536.

\item{[28]}
J.H.M.~Steenbrink and S.~Zucker, Variation of mixed Hodge structure,
I, Inv.\ Math. 80 (1985), 489--542.

\item{[29]}
R.~Thomas, Nodes and the Hodge conjecture, J. Alg.\ Geom.\ 14
(2005), 177--185.

\item{[30]}
G.N.~Tjurina, Locally semi-universal flat deformations of isolated
singularities of complex spaces, Izv.\ Akad.\ Nauk SSSR, Ser.\ Mat.\
33 (1969), 1026--1058.

\item{[31]}
C.~Voisin, Transcendental methods in the study of algebraic cycles,
in Lect. Notes in Math.\ vol.\ 1594, Springer, Berlin, 1994,
pp. 153--222.

{\sfont
\medskip
RIMS Kyoto University, Kyoto 606-8502 Japan

\vers}}\bye